\documentclass{SciPost}
\IfFileExists{mathdesign.sty}
  {\usepackage[bitstream-charter]{mathdesign}}
  {\usepackage[T1]{fontenc}\usepackage{lmodern,charter,amssymb}}
\usepackage{amsthm,booktabs,tikz}
\DeclareSymbolFont{usualmathcal}{OMS}{cmsy}{m}{n}
\DeclareSymbolFontAlphabet{\mathcal}{usualmathcal}
\hypersetup{colorlinks,linkcolor={scipostdeepblue},citecolor={blue!50!black},
  urlcolor={blue!80!black},pdftitle={A Systematic Construction of High-Corank Kac--Moody Algebras},
  pdfauthor={Simon Beaudoin, Quentin Bonnefoy, Alessio Marrani, Michel Rausch de Traubenberg and Victor Saulquin}}
\fancypagestyle{SPstyle}{\fancyhf{}
  \lhead{\colorbox{scipostblue}{\bfseries\color{white} ~SciPost Physics Proceedings }}
  \rhead{{\bfseries\color{scipostdeepblue} ~Submission }}
  \fancyfoot[C]{\textbf{\thepage}}}
\newtheorem{proposition}{Proposition}[section]
\newtheorem{corollary}[proposition]{Corollary}

\begin{document}
\pagestyle{SPstyle}
\begin{center}{\Large\bfseries\color{scipostdeepblue}
A Systematic Construction of High-Corank Kac--Moody Algebras\par}\end{center}
\begin{center}\textbf{Simon Beaudoin\textsuperscript{1$\star$},
Quentin Bonnefoy\textsuperscript{1$\dagger$},
Alessio Marrani\textsuperscript{2$\ddagger$},\\
Michel Rausch de Traubenberg\textsuperscript{1\S} and
Victor Saulquin\textsuperscript{1\P}}\end{center}
\begin{center}
{\bfseries 1} Universit\'e de Strasbourg, CNRS, IPHC UMR7178, F-67037 Strasbourg Cedex, France\\
{\bfseries 2} Centre for Mathematics and Theoretical Physics, University of Hertfordshire, Hatfield, UK\\[\baselineskip]
{\small $\star$ \href{mailto:simon.beaudoin@iphc.cnrs.fr}{simon.beaudoin@iphc.cnrs.fr}\quad
$\dagger$ \href{mailto:quentin.bonnefoy@iphc.cnrs.fr}{quentin.bonnefoy@iphc.cnrs.fr}\quad
$\ddagger$ \href{mailto:a.marrani@herts.ac.uk}{a.marrani@herts.ac.uk}\\
\S\ \href{mailto:michel.rausch@iphc.cnrs.fr}{michel.rausch@iphc.cnrs.fr}\qquad
\P\ \href{mailto:victor.saulquin@iphc.cnrs.fr}{victor.saulquin@iphc.cnrs.fr}}
\end{center}
\definecolor{palegray}{gray}{0.95}
\begin{center}
\colorbox{palegray}{%
  \begin{tabular}{@{}ll@{}}
 \begin{minipage}{0.20\textwidth}
    \includegraphics[width=28mm]{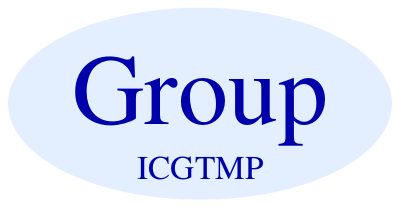}
  \end{minipage}
  &
  \begin{minipage}{0.752\textwidth}
    \begin{center}
    {\it 36th International Colloquium on Group Theoretical Methods in Physics}\\
    {\it Valladolid, 13--17 July 2026} \\
    \doi{10.21468/SciPostPhysProc.?}\\
    \end{center}
  \end{minipage}
\end{tabular}%
}
\end{center}
\section*{\color{scipostdeepblue}{Abstract}}
\textbf{\boldmath
We form a hexagonal prism from two six-cycles and iterate a doubling construction with a fixed positive integer edge multiplicity, producing a family of indefinite Kac--Moody algebras. The prism's spectrum determines the exact corank of every algebra in the family. Along each nonvanishing parity class, the corank grows exponentially, while its ratio to the matrix size tends to zero. For the twelve-node algebra we identify the two central directions, describe the first two graded levels, and characterize the integral central charges of its unitary highest-weight modules.}\par
\vspace{\baselineskip}

\begingroup
\setlength{\fboxsep}{2.4pt}
\setlength{\textwidth}{\dimexpr.975\linewidth-4\tabcolsep\relax}
\noindent\textcolor{white!90!black}{%
\fbox{\parbox{0.975\linewidth}{%
\textcolor{white!40!black}{\begin{tabular}{lr}%
  \begin{minipage}{0.6\textwidth}%
    {\small Copyright attribution to authors. \newline
    This work is a submission to SciPost Phys. Proc. \newline
    License information to appear upon publication. \newline
    Publication information to appear upon publication.}
  \end{minipage} & \begin{minipage}{0.4\textwidth}
    {\small Received Date \newline Accepted Date \newline Published Date}%
  \end{minipage}
\end{tabular}}
}}
}
\par\endgroup
\vspace{10pt}
\noindent\rule{\textwidth}{1pt}
{\setlength{\cftbeforesecskip}{.50em}\tableofcontents}
\noindent\rule{\textwidth}{1pt}
\vspace{10pt}

\section{Introduction}\label{sec:intro}
The corank of a generalized Cartan matrix $A$ can be read in two ways. On the Lie algebra side it is the dimension of the centre. When $A=(a_{ij})$ is symmetric, its corank is the
multiplicity of the adjacency eigenvalue $2$ of a loopless multigraph $G$. Give $G$ exactly $-a_{ij}$ edges between each pair of distinct vertices $i,j$, its adjacency matrix $\operatorname{Adj}(G)$ then satisfies
\begin{equation}\label{eq:graph}
A=2I-\operatorname{Adj}(G) \ ,\qquad
\dim(\ker (A))=\operatorname{mult}_{\operatorname{Adj}(G)}(2)\ .
\end{equation}
A connected finite-type diagram has corank zero, and an affine diagram has corank one. In indefinite type there is no upper bound independent of the number of vertices. Which connected diagrams give many central directions, and what becomes of the affine null root when the corank exceeds one?

Our starting point is the modified-cube construction of \cite{sqamv}. The first two diagrams are a three-vertex path and a six-cycle. Joining two copies of the cycle at corresponding vertices produces a hexagonal prism. Its Cartan matrix has order twelve, rank ten and one negative eigenvalue. After the prism, fix a positive integer $L$. At each step take two copies of the previous diagram and join each pair of corresponding vertices by $L$ parallel edges. Each old eigenvalue splits into two translates. Their multiplicities give an exact expression for the corank and explain why a change from one edge to two changes the behaviour at successive steps. Reference~\cite{sqamv} also treats a second family built from weighted cubes. For the ordinary even-dimensional cubes in that family, the same binomial corank formula was recorded by Terence C in January 2025 \cite{TerenceC2025}.

The prism also gives a small enough example to follow the two central directions through the Lie algebra. In Section~\ref{sec:algebra} we mark two adjacent vertices. Their simple-root coefficients give a bigrading, and their neighbours determine the extremal weights of the first graded pieces. At the
second level, the edge between the marked vertices decides which Serre relations can occur. This gives a direct calculation of the mixed piece. We then ask which pairs of central charges occur in unitary highest-weight modules and exhibit infinitely many dominant weights above each integral pair.

One motivation for studying indefinite diagrams comes from the use of Kac--Moody symmetries in gravity and supergravity, beginning with Julia's work \cite{Julia85}. The $E_{10}$ correspondence compares low levels of a coset model with the equations of eleven-dimensional supergravity near a spacelike singularity \cite{DamourHenneauxNicolai2002,DamourHenneauxNicolai2003}.
West proposed $E_{11}$ as a symmetry underlying M-theory \cite{WestE11}. Here, the physical parallel concerns the use of a level decomposition. The existence of a field-theoretic realization for the diagrams studied below is open.

\section{Cartan matrices, centres and null vectors}\label{sec:gcm}

Let $A=(a_{ij})$ be an indecomposable symmetrizable generalized Cartan matrix of order $r$. Thus, $a_{ii}=2$, $a_{ij}\in\mathbb Z_{\leq0}$ for $i\ne j$, and $a_{ij}=0$ precisely when $a_{ji}=0$. Indecomposability means that the associated diagram is connected. Fix a positive diagonal matrix $D$ such that $DA$ is
symmetric. The type is finite if $DA$ is positive definite and affine if it is positive semidefinite and singular. In the latter case its kernel is one-dimensional. In all other cases, note that the type is indefinite \cite{Kac,Moody,Kac2,moL}.

Put $s = \dim(\ker(A))$. A minimal realization consists of a space $\mathfrak h$ of dimension $r+s$, independent simple coroots $h_i\in\mathfrak h$ and independent simple roots $\alpha_j\in\mathfrak h^*$ with $\alpha_j(h_i)=a_{ij}$ for $1\leq i,j\leq r$. Write $\mathfrak h'=\bigoplus_{i=1}^r\mathbb C h_i$ for the coroot span. The complex algebra $\mathfrak g(A)$ is generated by $\mathfrak h$ and $e_i,f_i$ subject to
\begin{align}
 [h,h']&=0\ , & [e_i,f_j]&=\delta_{ij}h_i\ ,\notag\\
 [h,e_i]&=\alpha_i(h)e_i\ , & [h,f_i]&=-\alpha_i(h)f_i\ ,\label{eq:chevalley} \\
 (\operatorname{ad}(e_i))^{1-a_{ij}}e_j&=0\ , &
 (\operatorname{ad}(f_i))^{1-a_{ij}}f_j&=0\qquad(i\ne j)\ .\notag
\end{align}
Here $h,h'\in\mathfrak h$. For symmetrizable $A$, these relations give the Chevalley--Serre presentation \cite{Kac2,GabberKac}. Accounts of the
construction and its root system can also be found in \cite{Mdo,Mar}. The centre of $\mathfrak g(A)$ consists of the elements $K_u=\sum_{i=1}^r u_i h_i\in\mathfrak h'$ with $A^{\mathsf T}u=0$. On the root lattice $Q=\bigoplus_i\mathbb Z\alpha_i$, we denote the invariant symmetric bilinear form by $\langle\,\cdot\,,\cdot\,\rangle$ and normalize it by $\langle\alpha_i,\alpha_j\rangle=(DA)_{ij}$. Its radical in $Q\otimes\mathbb R$ consists of the $\delta_v=\sum_i v_i\alpha_i$ with $Av=0$. Thus the centre uses the left kernel and the root radical uses the right kernel. For the symmetric matrices of \eqref{eq:graph} the same coefficient vector gives both.

\begin{proposition}\label{prop:radical}
For an indecomposable symmetrizable GCM, the existence of a nonzero real vector $v\geq0$ in its kernel forces affine type. For an indecomposable indefinite GCM, every nonzero real kernel vector has coefficients of both signs. Consequently the integral radical contains no root in indefinite type.
\end{proposition}
\begin{proof}
$N=2I-A$ is an irreducible nonnegative matrix and $Nv=2v$ since $v\in\ker(A)$. If one coordinate of $v\geq0$ vanishes, the same holds for each adjacent coordinate, hence for all of them. Thus $v>0$ and Perron--Frobenius gives $\rho(N)=2$ as a simple
eigenvalue. In the symmetrizable case, $N$ is similar to a real symmetric matrix. It follows that $D^{1/2}AD^{-1/2}$, and hence its congruent matrix $DA$, is positive semidefinite of corank one. This is the affine criterion \cite{mo,Kac2}. The coefficient vector of any root is either nonnegative or nonpositive in the simple-root basis. The assertion about the integral radical follows.
\end{proof}

A vector can have square zero without belonging to the radical. The distinction already appears in the prism, where the six-cycle carries an isotropic imaginary root that pairs nontrivially with the other layer (Section~\ref{sec:realization}).

\section{Spectra of the modified cubes}\label{sec:spectra}

\subsection{The seed and the doubling operation}\label{sec:seed}

\begin{figure}[htbp]
\centering
\begin{tikzpicture}[font=\small,vertex/.style={circle,draw=scipostdeepblue,fill=white,
 inner sep=0pt,minimum size=4.6mm},marked/.style={vertex,fill=scipostdeepblue,text=white},
 edge/.style={draw=scipostdeepblue,line width=.8pt},rung/.style={draw=black!58,dashed,line width=.8pt}]
 \begin{scope}[xshift=0cm]
 \foreach \i/\x in {1/0,2/.85,3/1.7}{\coordinate (p\i) at (\x,0);}
 \draw[edge] (p1)--(p2)--(p3);
 \foreach \i in {1,...,3}{\node[vertex] at (p\i){\i};}
 \node at (.85,-2.23) {$G_1$};
 \end{scope}
 \begin{scope}[xshift=4.5cm]
 \foreach \i in {1,...,6}{\coordinate (u\i) at ({90-60*(\i-1)}:1);}
 \draw[edge] (u1)--(u2)--(u3)--(u4)--(u5)--(u6)--cycle;
 \foreach \i in {1,...,6}{\node[vertex] at (u\i){\i};}
 \node at (0,-2.23) {$G_2=C_6$};
 \end{scope}
 \begin{scope}[xshift=8.3cm]
 \foreach \i in {1,...,6}{
   \coordinate (o\i) at ({90-60*(\i-1)}:1.65);
   \coordinate (i\i) at ({90-60*(\i-1)}:.82);
   \draw[rung] (o\i)--(i\i);}
 \draw[edge] (o1)--(o2)--(o3)--(o4)--(o5)--(o6)--cycle;
 \draw[edge] (i1)--(i2)--(i3)--(i4)--(i5)--(i6)--cycle;
 \foreach \i in {1,...,6}{\node[vertex] at (o\i){\i};}
 \foreach \i/\j in {1/7,2/8,3/9,4/10}{\node[vertex] at (i\i){\j};}
 \node[marked] at (i5){11}; \node[marked] at (i6){12};
 \node at (0,-2.23) {$G_3=C_6\mathbin{\square}K_2$};
 \end{scope}
\end{tikzpicture}
\caption{The three seed diagrams. The dashed edges of the prism join $i$ to $i+6$ and have multiplicity one. The filled vertices are those used for the bigrading in Section~\ref{sec:grading}.}\label{fig:seed}
\end{figure}
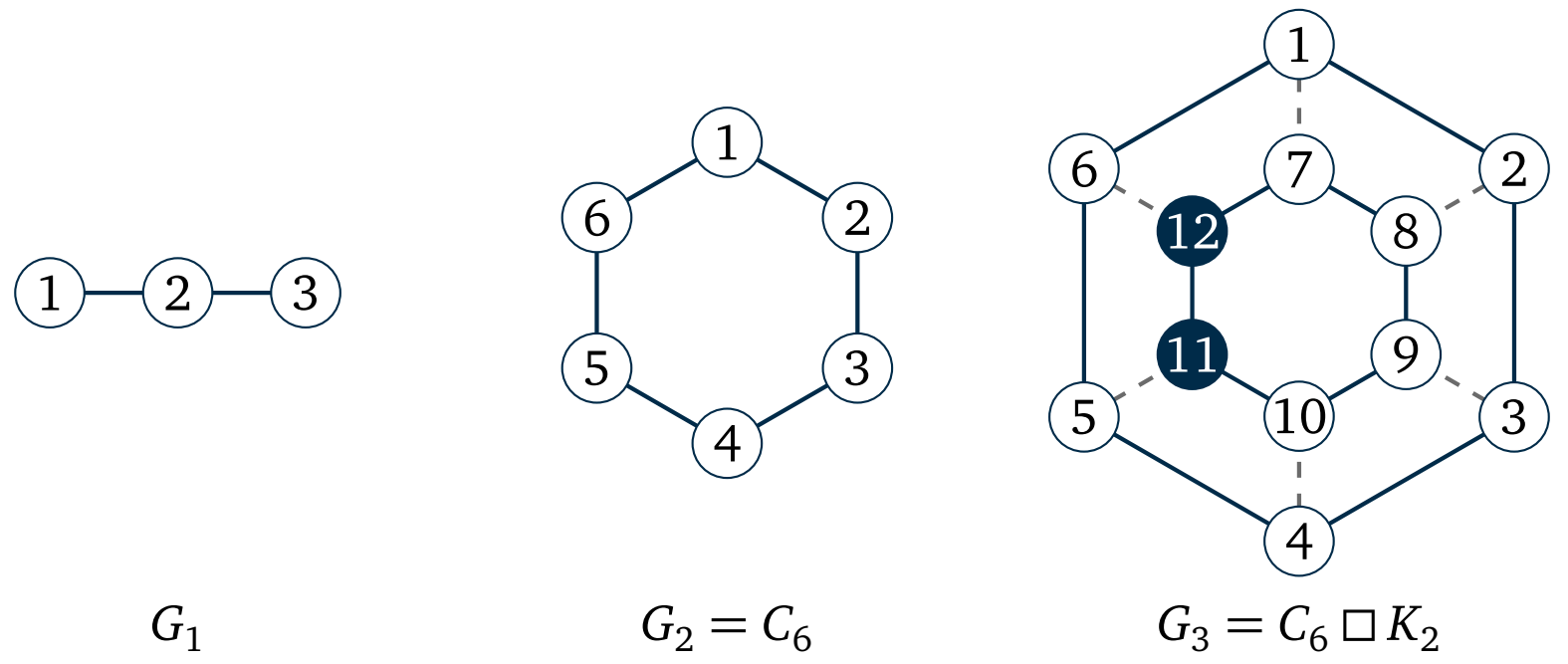

Write $C_6$ for the six-cycle and $K_2^{(L)}$ for the two-vertex multigraph with $L\geq1$ edges between its vertices. We write $K_2=K_2^{(1)}$. Deleting one vertex of a square gives the three-vertex path $G_1$. In the next step two copies of this path are joined at opposite ends to form $G_2=C_6$, also obtained by deleting two antipodal vertices of a cube. The path gives the finite algebra $\mathfrak a_3\simeq\mathfrak{sl}_4$,
and the six-cycle has untwisted affine type $\widehat{\mathfrak a}_5$ \cite{Ca}. Throughout the paper, a subscript on a \emph{matrix} $A_n$ denotes the construction step. At step three, join corresponding nodes of two six-cycles by one edge. This gives the prism $G_3=C_6\mathbin{\square}K_2$, with the numbering in
Figure~\ref{fig:seed}. Each square face of the prism has affine type
$\widehat{\mathfrak a}_3$. The parameter $L$ enters at the first doubling after the prism.

For the Cartesian product $G\mathbin{\square}H$, the adjacency matrix is the sum of $\operatorname{Adj}(G)\otimes I$ and $I\otimes\operatorname{Adj}(H)$
\cite{BrouwerHaemers}. If $A$ is the Cartan matrix of $G$, the Cartan matrix
of $G\mathbin{\square}K_2^{(L)}$, in layer order, is
\begin{equation}\label{eq:double}
 \begin{pmatrix}A&-LI\\-LI&A\end{pmatrix}\ .
\end{equation}
For $Ax=\mu x$, the vectors $(x,x)$ and $(x,-x)$ have eigenvalues $\mu-L$ and $\mu+L$, respectively. They account for the entire spectrum. Thus if $P(z)=\det(zI-A)$, the characteristic polynomial after one doubling is $P(z-L)P(z+L)$. The same argument works for any initial connected loopless multigraph. The prism construction fixes the initial multiplicities explicitly.

For the cycle, $A_2=2I_6-\operatorname{Adj}(C_6)$ has eigenvalues
$2-2\cos(2\pi j/6)$ for $0\leq j<6$. Shifting these by $\pm1$ gives
\begin{equation}\label{eq:seed}
 A_3=\begin{pmatrix}A_2&-I_6\\-I_6&A_2\end{pmatrix}\ ,\qquad
 \begin{array}{c|rrrrrrr}
 \mu&-1&0&1&2&3&4&5\\\hline
 d_3(\mu)&1&2&1&4&1&2&1
 \end{array}
\end{equation}
where $d_3(\mu)$ is the multiplicity. Thus $A_3$ has inertia
$(9,1,2)$, since it has $9$ positive eigenvalues with multiplicity, $1$ negative eigenvalue and $2$ zero eigenvalues. The graph at step $n=q+3$ is
\begin{equation}\label{eq:product}
G_n^{(L)}=C_6\mathbin{\square}K_2\mathbin{\square}
(K_2^{(L)})^{\square q}\ ,\qquad |V(G_n^{(L)})|=12 \cdot 2^q
\end{equation}
where $V(G)$ denotes the vertex set of $G$. We write
$A_n^{(L)}=2I-\operatorname{Adj}(G_n^{(L)})$ for its Cartan matrix.
In particular $A_3^{(L)}=A_3$ for every $L$. The empty product at $q=0$ yields the prism. To see the modified cube
directly, label the vertices of the ordinary $(n+1)$-cube by
$\{0,1\}^3\times\{0,1\}^{n-2}$. Delete the vertices with first three
coordinates $000$ or $111$. The induced graph is the underlying simple
graph in \eqref{eq:product}. Of its last $n-2$ coordinate directions,
one carries single edges and the other $q$ carry multiplicity $L$.
Every graph in \eqref{eq:product}
is connected and regular of degree $3+qL$, with edges counted with
multiplicity. Its constant vector is an eigenvector of $A_n^{(L)}$ of
eigenvalue $-1-qL$. Thus all these matrices are indefinite, including the
invertible members of the family.

\subsection{Exact multiplicities and parity}\label{sec:parity}

For $q\geq0$, set
\begin{equation}\label{eq:binomial}
b_q(t)=\begin{cases}
\displaystyle\binom{q}{(q+t)/2}\ , & (q+t)/2\in\{0,1,\ldots,q\}\ ,\\[2pt]
0\ ,&\text{otherwise}.
\end{cases}
\end{equation}
In particular, $b_q(t)=b_q(-t)$. This convention allows to give a single compact formula that, with $d_3(\mu)$, characterizes the whole spectrum. It embeds parity and range conditions, since it vanishes when $t \not\equiv q$ mod $2$ or $|t| > q$.

\begin{proposition}[{\cite{sqamv}}]\label{prop:count}
Let $S=\{-1,0,1,2,3,4,5\}$ be the spectrum of $A_3$. At step $n=q+3$, the spectral multiplicity of $\lambda$ and the corank $c_n(L)$ are
\begin{align}
 d_n^{(L)}(\lambda)&=\sum_{\mu\in S} d_3(\mu)\,
                 b_q\big((\lambda-\mu)/L\big)\ ,\label{eq:count}\\
 c_n(L)&=\sum_{\mu\in S}d_3(\mu)\,b_q(\mu/L)\ .\label{eq:corank}
\end{align}
\end{proposition}
\begin{proof}
Each tensor product of a seed eigenvector and $q$ symmetric or antisymmetric two-component vectors is an eigenvector of the product. Choosing $j$ of the $q$ shifts to be $+L$ gives $\lambda=\mu+L(2j-q)$ with multiplicity $d_3(\mu)\binom qj$. Summing contributions at the same $\lambda$ proves \eqref{eq:count}. Setting $\lambda=0$ in \eqref{eq:count} and using
$b_q(-t)=b_q(t)$ gives \eqref{eq:corank}.
\end{proof}
Equivalently, all multiplicities are Laurent coefficients of
\begin{equation}\label{eq:laurent}
\sum_{\lambda}d_n^{(L)}(\lambda)x^\lambda
 =\big(x^{-1}+2+x+4x^2+x^3+2x^4+x^5\big)(x^{-L}+x^L)^q \ .
\end{equation}
For instance, at $n=4$, one reads $c_4(1)=d_3(-1)+d_3(1)=2$ and
$c_4(2)=d_3(-2)+d_3(2)=4$, taking $d_3$ to be zero off $S$.

For the corank, it suffices to keep the seed eigenvalues divisible by $L$. Table~\ref{tab:classes} separates the resulting five possibilities. For $L=3$ and $5$ the only nonzero seed eigenvalue that survives is $L$ itself, with multiplicity one. Hence the two sequences agree term by term. For $L\geq6$ only the zero eigenvalue survives. In particular
\begin{equation}\label{eq:largeL}
 c_{q+3}(L)=
 \begin{cases}2\binom{q}{q/2}\ ,&q\text{ even}\ ,\\0\ ,&q\text{ odd}\ ,\end{cases}
 \qquad L\geq6 \ .
\end{equation}
This also explains why a connected indefinite diagram can alternate between
positive corank and full rank as more vertices are added.

\begin{table}[htbp]
\centering
\caption{Exact corank and asymptotic coefficients for every $L$. In the second column $b_j$ abbreviates $b_q(j)$ with $q=n-3$.
The last two columns give the coefficients of Corollary~\ref{cor:growth} for even and odd $q$, respectively.}\label{tab:classes}
\begin{tabular}{@{}llrr@{}}
\toprule
$L$ & $c_{q+3}(L)$ & $M_{L,0}$ & $M_{L,1}$\\
\midrule
$1$ & $2b_0+2b_1+4b_2+b_3+2b_4+b_5$ & 8 & 4\\
$2$ & $2b_0+4b_1+2b_2$                 & 4 & 4\\
$3,5$ & $2b_0+b_1$                       & 2 & 1\\
$4$ & $2b_0+2b_1$                        & 2 & 2\\
$\geq6$ & $2b_0$                          & 2 & 0\\
\bottomrule
\end{tabular}
\end{table}

\begin{figure}[htbp]
\centering
\includegraphics[width=.96\linewidth]{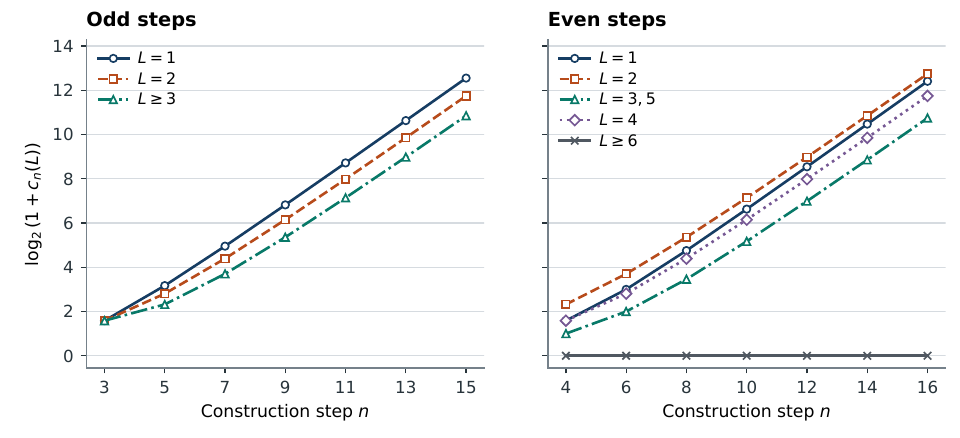}
\caption{Coranks of the modified cubes on odd and even construction steps. At each step $n$, the marker has ordinate $\log_2(1+c_n(L))$, where $c_n(L)$ is the exact integer corank given by \eqref{eq:corank}. A zero corank is plotted at zero. The curves for all $L\geq3$ coincide on odd steps.\label{fig:coranks}}
\end{figure}

Some differences between successive steps are easier to see before taking a limit. Pascal's identity applied to the first two rows of
Table~\ref{tab:classes} gives
\begin{align}
c_n(1)&=
 \begin{cases}
  \binom{n-2}{(n-2)/2}+\binom{n-2}{(n-6)/2}\ ,&n\text{ even}\ ,\\
  2\binom{n-1}{(n-3)/2}\ ,&n\text{ odd}\ ,
 \end{cases}\label{eq:one}\\
c_n(2)&=
 \begin{cases}
  2\binom{n-2}{(n-2)/2}\ ,&n\text{ even}\ ,\\
  \binom{n-1}{(n-1)/2}\ ,&n\text{ odd}\ .
 \end{cases}\label{eq:two}
\end{align}
These formulas hold for $n\geq3$, with a binomial coefficient taken as zero when its lower index is outside its range. For $m\geq2$, they yield
\begin{equation}\label{eq:ratios}
\begin{array}{c|cc}
 & c_{2m+1}(L)/c_{2m}(L) & c_{2m+2}(L)/c_{2m+1}(L)\\[2pt]\hline
 L=1 & \dfrac{2m(2m-1)}{m^2-m+1} & \dfrac{m^2+m+1}{m(m+2)}\\[7pt]
 L=2 & \dfrac{2m-1}{m} & 2
\end{array}
\end{equation}
Thus, for $L=1$, the corank falls at every transition from $2m+1$ to $2m+2$ with $m\geq2$ and rises on the intervening steps. For $L=2$, it increases at every step starting at the prism. Table~\ref{tab:values} gives the integer values, and Figure~\ref{fig:coranks} shows the two parity sequences for every edge multiplicity.

\begin{table}[htbp]
\centering
\caption{The first twelve steps. Steps one and two precede the introduction of $L$. From step three onwards the order is $12\cdot2^{n-3}$.}\label{tab:values}
\begin{tabular}{@{}r r r r r r r@{}}\toprule
$n$ & Order & $L=1$ & $L=2$ & $L=3,5$ & $L=4$ & $L\geq6$\\\midrule
1 & 3    & 0 & 0 & 0 & 0 & 0\\
2 & 6    & 1 & 1 & 1 & 1 & 1\\
3 & 12   & 2 & 2 & 2 & 2 & 2\\
4 & 24   & 2 & 4 & 1 & 2 & 0\\
5 & 48   & 8 & 6 & 4 & 4 & 4\\
6 & 96   & 7 & 12 & 3 & 6 & 0\\
7 & 192  & 30 & 20 & 12 & 12 & 12\\
8 & 384  & 26 & 40 & 10 & 20 & 0\\
9 & 768  & 112 & 70 & 40 & 40 & 40\\
10 & 1536 & 98 & 140 & 35 & 70 & 0\\
11 & 3072 & 420 & 252 & 140 & 140 & 140\\
12 & 6144 & 372 & 504 & 126 & 252 & 0\\
\bottomrule\end{tabular}
\end{table}

\begin{corollary}\label{cor:growth}
Fix $L$. On the parity class $q\equiv p\pmod2$, let
$M_{L,p}=\sum_{\mu\in S,\ \mu/L\in 2\mathbb Z+p}d_3(\mu)$.
If $M_{L,p}>0$, then, along that class,
\begin{equation}\label{eq:asymptotics}
 c_{q+3}(L)=M_{L,p}\,2^q\sqrt{\frac{2}{\pi q}}\,
 (1+O(q^{-1}))\ .
\end{equation}
If $M_{L,p}=0$, the corank vanishes throughout the class. In particular, $c_n(L)/|V(G_n^{(L)})|\to0$ for every $L$.
\end{corollary}
\begin{proof}
Fix an integer $t$ of the same parity as $q$. Stirling's formula gives
\begin{equation}
 b_q(t)=2^q\sqrt{2/(\pi q)}\,\bigl(1+O(q^{-1})\bigr)\ .
\end{equation}
All terms of the opposite parity or with nonintegral $t$ in
\eqref{eq:corank} vanish. There are only seven seed eigenvalues, so the estimate can be summed term by term. Division by $12\cdot 2^q$ gives the last
assertion.
\end{proof}
The factor $q^{-1/2}$ has a simple interpretation. There are $2^q$ sign choices in the doubling construction. A kernel vector can only occur when the sum of the $q$ signs equals one of finitely many values fixed by the seed eigenvalues. Any single value of this sum occurs for only about a $q^{-1/2}$ proportion of the $2^q$ sign choices.

\section{The twelve-node algebra}\label{sec:algebra}

\subsection{A minimal realization}\label{sec:realization}

We use the numbering of Figure~\ref{fig:seed}, in which
vertices 11 and 12 are marked. The kernel of $A_3$ is
the symmetric lift $(u,u)$ of the eigenspace of $A_2$ at eigenvalue one. Indeed, the antisymmetric lift would require an eigenvalue $-1$ of $A_2$.
Define two column vectors by
\begin{align}
 v_1&=(u_1,u_1)\ ,&u_1&=(-1,-1,0,1,1,0)\ ,\notag\\
 v_2&=(u_2,u_2)\ ,&u_2&=(1,0,-1,-1,0,1)\ .\label{eq:v}
\end{align}
Both satisfy $A_2u_a=u_a$, with the six-tuples here read as columns. The last two rows of the matrix $(v_1\ v_2)$ form $I_2$, so $v_1,v_2$ are independent and span the two-dimensional kernel. The coefficients of a kernel vector in this basis are its last two coordinates. Thus $v_1,v_2$ also form a $\mathbb Z$-basis of $\ker(A_3)\cap\mathbb Z^{12}$. With this choice, the central elements and the radical vectors have matching coordinates,
\begin{equation}\label{eq:central}
K_a=\sum_{i=1}^{12}(v_a)_i h_i\ ,\qquad
\delta_a=\sum_{i=1}^{12}(v_a)_i\alpha_i\qquad (a=1,2)\ .\end{equation}
The $K_a$ span the centre of $\mathfrak g(A_3)$ and the $\delta_a$ span the root radical. Every nonzero integral combination of the latter has mixed signs by Proposition~\ref{prop:radical}. The sum
$\eta=\alpha_1+\cdots+\alpha_6$ illustrates the other possibility for a null vector. The subdiagram on nodes $1,\ldots,6$ has affine type $\widehat{\mathfrak a}_5$, and the six Chevalley triples generate its derived affine algebra. The primitive positive imaginary root of this affine root system is $\eta$, hence $\eta$ is a root of $\mathfrak g(A_3)$. Let $c_\eta=\binom{\mathbf1_6}{0_6}$ be its coefficient column in the simple-root
basis, with $\mathbf1_6=(1,\ldots,1)^{\mathsf T}$.
Equation~\eqref{eq:seed} gives $A_3c_\eta=\binom{0_6}{-\mathbf1_6}$. Thus
\begin{equation}
\langle\eta,\eta\rangle=c_\eta^{\mathsf T}A_3c_\eta=0\ ,\qquad
\langle\eta,\alpha_{6+j}\rangle=-1\quad(1\leq j\leq6)\ .
\end{equation}
Each rung pairs $\eta$ with a vertex of the opposite layer, so this isotropic root lies outside the radical of the full root lattice.

Let $H$ be the principal submatrix on nodes $1,\ldots,10$, and write
\begin{equation}\label{eq:gram}
A_3=\begin{pmatrix}H&B\\ B^{\mathsf T}&C\end{pmatrix}\ ,\qquad
\mathcal C=\begin{pmatrix}H&B&0\\ B^{\mathsf T}&C&I_2\\0&I_2&0\end{pmatrix}\ .
\end{equation}
Here, $C=\left(\begin{smallmatrix}2&-1\\-1&2\end{smallmatrix}\right)$. If $\epsilon_i$ denotes the $i$th coordinate vector in $\mathbb R^{10}$, the columns of $B$ are $-(\epsilon_5+\epsilon_{10})$ and
$-(\epsilon_6+\epsilon_7)$. Direct evaluation gives $\det H=-80$.
The coordinates of \eqref{eq:v} give another way to establish nondegeneracy. The span of the first ten simple roots is a complement to the radical, and the form restricted to a complement of its radical is nondegenerate. Its inertia is $(9,1,0)$ by \eqref{eq:seed}.

Adjoin $d_1,d_2$ to the coroot span and set
$\alpha_i(d_a)=\delta_{i,10+a}$. Order the basis with all twelve $h_i$
first, followed by $d_1,d_2$. In this basis, the invariant form on the Cartan subalgebra can be chosen to have Gram matrix $\mathcal C$. With $K_1,K_2$ in place of $h_{11},h_{12}$ this becomes
\begin{equation}\label{eq:split}
H\ \oplus\ \begin{pmatrix}0&I_2\\ I_2&0\end{pmatrix}\ .
\end{equation}
In particular, $\langle K_a,d_b\rangle=\delta_{ab}$, and $\delta_a(d_b)=\delta_{ab}$. The extended form has inertia $(11,3,0)$ on the real span of this basis. The twelve simple roots are independent on the resulting fourteen-dimensional Cartan subalgebra. This is a minimal realization. It also fixes the signs
of the grading operators,
\begin{equation}\label{eq:degree}
 [d_a,e_i]=\delta_{i,10+a}e_i\ ,\qquad
 [d_a,f_i]=-\delta_{i,10+a}f_i\ .
\end{equation}
The $K_a$ commute with the whole algebra. The degree operators
$d_a$ belong to $\mathfrak g(A_3)$ and act by derivations on its derived algebra $\mathfrak g'(A_3)$. Indeed, $\mathfrak h\cap\mathfrak g'(A_3)=\mathfrak h'$, and each $\delta_b$ vanishes on $\mathfrak h'$ while $\delta_b(d_a)=\delta_{ab}$. Thus the $d_a$ lie outside the derived algebra.

\subsection{The first two levels of the bigrading}\label{sec:grading}

For a root $\beta=\sum_i t_i\alpha_i$, set
$\deg\beta=(t_{11},t_{12})$ and put $\mathfrak h$ in degree zero. Write $\mathfrak g_{u,v}$ for the corresponding summands of $\mathfrak g(A_3)$. Root spaces occur only where $u,v$ are both nonnegative or both nonpositive. The algebra $\mathfrak l=\mathfrak g(H)$ on the ten unmarked nodes is an indefinite Kac--Moody algebra with nonsingular Cartan matrix. The zero component is
\begin{equation}\label{eq:zero}
\mathfrak g_{0,0}=\mathfrak l\oplus\mathbb CK_1\oplus\mathbb CK_2
                      \oplus\mathbb Cd_1\oplus\mathbb Cd_2\ .
\end{equation}
In fact, a positive root of degree zero has support entirely on the unmarked nodes, and the same reasoning applies to its negative. The remaining Cartan directions commute with $\mathfrak l$.

Let $p_a$ denote the first ten entries of $v_a$, and let $\omega_i$ denote the fundamental weights of $\mathfrak l$, with $\omega_i(h_j)=\delta_{ij}$ for $i,j\leq10$. Identify these weights with elements in the span of the first ten simple roots using $H^{-1}$. Define $J_1=\{5,10\}$ and $J_2=\{6,7\}$. Write $B_{\bullet a}$ for the $a$th column of $B$ in \eqref{eq:gram}.
Since $Hp_a=-B_{\bullet a}$, the weights
$\psi_a=\sum_{i=1}^{10}(p_a)_i\alpha_i$ are
\begin{equation}\label{eq:signs}
 \psi_1=\omega_5+\omega_{10}\ ,\qquad
 \psi_2=\omega_6+\omega_7\ ,\qquad
 \alpha_{10+a}=\delta_a-\psi_a \ .
\end{equation}
Thus, the $\mathfrak l$-weight of $x_a=e_{10+a}$ is $-\psi_a$.
For example $[h_5,e_{11}]=-e_{11}$ although $\psi_1(h_5)=1$.
This sign follows directly from the edge $5$--$11$.

Denote the irreducible integrable lowest-weight module of extremal weight
$-\psi_a$ by $V_a$ and its restricted dual highest-weight module by
$V_a^\vee$ (the direct sum of the dual weight spaces with the contragredient
action). Then, as $\mathfrak l$-modules,
\begin{equation}\label{eq:firstlevel}
\mathfrak g_{1,0}\simeq V_1\ ,\qquad \mathfrak g_{0,1}\simeq V_2\ ,\qquad
 \mathfrak g_{-1,0}\simeq V_1^\vee\ ,\quad
 \mathfrak g_{0,-1}\simeq V_2^\vee \ .
\end{equation}
Each positive component in \eqref{eq:firstlevel} is generated, as an $\mathfrak l$-module, by its marked simple generator. A positive bracket word of bidegree $(1,0)$ contains $e_{11}$ exactly once. Repeated use of the Jacobi identity writes it as a sum of iterated brackets of unmarked generators with
$e_{11}$. The same argument at $(0,1)$ starts from $e_{12}$. In particular,
\begin{equation}
\mathfrak g_{1,0}=U(\mathfrak n_+(\mathfrak l))\cdot e_{11}\ ,\qquad
 \mathfrak g_{0,1}=U(\mathfrak n_+(\mathfrak l))\cdot e_{12}\ .
\end{equation}
Every $f_j$ with $j\leq10$ kills $x_a$, and the Serre relations give $(\operatorname{ad}e_j)^{1-(A_3)_{j,10+a}}x_a=0$. Together with the weight $-\psi_a$, these are the defining lowest-weight relations for $V_a$ \cite[Chaps.~9--10]{Kac2}. They give a surjection from $V_a$ onto the component generated by $x_a$. Since $x_a\ne0$ and $V_a$ is irreducible, this map is an isomorphism. The negative pieces follow by the Chevalley involution, or by the
nondegenerate invariant pairing between opposite root spaces.

The first six nodes distinguish $V_a$ from its dual. Their coroot sum $h_{\mathrm{cyc}}=\sum_{i=1}^6 h_i$ acts on $x_a$ by $-1$.
We have $\alpha_j(h_{\mathrm{cyc}})=0$ for $j\leq6$ and
$\alpha_j(h_{\mathrm{cyc}})=-1$ for $7\leq j\leq10$. Raising by unmarked generators therefore keeps the spectrum of $h_{\mathrm{cyc}}$ on $V_a$ in $\mathbb Z_{\leq-1}$. On the restricted dual, all eigenvalues are positive, so the two modules are nonisomorphic.

The same presentation determines what happens at total level two. Use $e_jx_a$ to mean the adjoint action $[e_j,x_a]$ in $V_a$, and let $U(\mathfrak l)$ act diagonally on tensors. The following description specifies the relation submodules as well as the ambient tensor products.

\begin{proposition}\label{prop:second}
As $\mathfrak l$-modules, the Lie bracket induces the identifications
\begin{equation}
\mathfrak g_{1,1} \simeq V_1\otimes V_2\ ,\label{eq:mixed}
\end{equation}
\begin{equation}
\mathfrak g_{2,0}\simeq \bigwedge\nolimits^2 V_1/R_1\ , \qquad 
\mathfrak g_{0,2}\simeq \bigwedge\nolimits^2 V_2/R_2\ ,\label{eq:pure}
\end{equation}
\begin{equation}
R_a=\sum_{j\in J_a}U(\mathfrak l) \bigl(x_a\wedge(e_jx_a)\bigr)\ .
\end{equation}
The displayed generators of $R_a$ are lowest-weight vectors of weights $-2\psi_a+\alpha_j|_{\mathfrak h_{\mathfrak l}}$.
\end{proposition}
\begin{proof}
Start with the positive Serre presentation and impose first the relations among unmarked generators. Write $\mathfrak n=\mathfrak n_+(\mathfrak l)$ for the resulting unmarked algebra. The positive algebra at this stage is $\mathfrak n\ltimes F(W)$, where $F$ denotes the free Lie algebra and
$W=U(\mathfrak n)\otimes(\mathbb Cx_1\oplus\mathbb Cx_2)$. Indeed, the left action of $\mathfrak n$ on $W$ extends by derivations to $F(W)$. Conversely, any images of the generators satisfying the unmarked relations send $u\otimes x_a$ to the iterated adjoint action of $u$ on the image of $x_a$. The universal property of $F(W)$ then gives the inverse presentation. The relations with just one occurrence of a marked generator are linear in $W$. Factoring out their $\mathfrak n$-span replaces $F(W)$ by
$F(V_1\oplus V_2)$, by \eqref{eq:firstlevel}. The remaining Serre relations generate a homogeneous ideal there, together with their $\mathfrak n$-translates.

In total degree two, the free components are $\bigwedge^2 V_1$, $V_1\otimes V_2$ and $\bigwedge^2 V_2$. A defining
relation with the same marked generator twice has the form
$[x_a,[x_a,e_j]]=0$ for $j\in J_a$. Its image in the exterior square, up to sign, is $x_a\wedge(e_jx_a)$. All its unmarked positive translates must be factored out. Each generator of this relation space is killed by every $f_i$ with $i\leq10$, since $f_i(e_jx_a)=-\delta_{ij}h_jx_a$ is a multiple of $x_a$. Its positive translates therefore span its full $\mathfrak l$-submodule. Finally, the two marked vertices are adjacent. Their mutual Serre relations have bidegrees $(2,1)$ and $(1,2)$. In the
positive presentation, bracketing can only increase a marked coefficient, so these relations contribute nothing in bidegree $(1,1)$.
\end{proof}
The two summands defining each $R_a$ can themselves be identified. Since $j\in J_a$, the root string through $x_a$ gives $e_jx_a\ne0$ and the exterior product in the proposition is nonzero. Its lowest weight is $-\gamma_{a,j}$, where $\gamma_{a,j}=2\psi_a-\alpha_j|_{\mathfrak h_{\mathfrak l}}$ is dominant. Each vector generates an irreducible integrable lowest-weight module. For fixed $a$, the two weights differ, since $H$ is nonsingular, so the sum is direct. In fundamental-weight coordinates the two pairs are
\begin{align}
 (\gamma_{1,5},\gamma_{1,10})
  &=(\omega_4+\omega_6+2\omega_{10},\ \omega_4+2\omega_5+\omega_9)\ ,\notag\\
 (\gamma_{2,6},\gamma_{2,7})
  &=(\omega_1+\omega_5+2\omega_7,\ \omega_1+2\omega_6+\omega_8)\ .
 \label{eq:relationsweights}
\end{align}
In the mixed component, $[e_{11},e_{12}]$ is the image of $x_1\otimes x_2$ and generates the lowest-weight constituent of weight $-\psi_1-\psi_2$. The negative components have the corresponding contragredient descriptions.

\subsection{Unitary highest weights and central characters}\label{sec:weights}

We use the conjugate-linear anti-involution defined by $e_i^*=f_i$,
$h_i^*=h_i$ and $d_a^*=d_a$. The compact real form consists of the elements $X$ satisfying $X^*=-X$. By a unitary highest-weight module, we mean an algebraic highest-weight module carrying a positive-definite Hermitian form for which $X^*$ acts as the adjoint of $X$. Let $L(\Lambda)$ be the \emph{irreducible} module
of highest weight $\Lambda$, and $v_\Lambda$ a unit highest-weight vector. Put $\ell_i=\Lambda(h_i)$ for all twelve Dynkin labels and use the coordinates
\begin{equation}\label{eq:weightcoords}
n_i=\ell_i\ (1\leq i\leq10)\ ,\qquad
k_a=\Lambda(K_a)\ ,\qquad m_a=\Lambda(d_a)\quad(a=1,2)\ .\end{equation}
Equation~\eqref{eq:central} gives
$h_{10+a}=K_a-\sum_{i=1}^{10}(p_a)_i h_i$. The two missing Dynkin labels are therefore $\ell_{10+a}=k_a-p_a^{\mathsf T}n$, where $n=(n_1,\ldots,n_{10})^{\mathsf T}$. Each $K_a$ acts by the scalar $k_a$ throughout the module. On a lowering monomial applied to $v_\Lambda$, the $d_a$-eigenvalue is $m_a-q_a$, where $q_a$ counts the occurrences of $f_{10+a}$.

\begin{proposition}\label{prop:unitary}
$L(\Lambda)$ is unitary for the compact form if and only if
\begin{equation}\label{eq:unitary}
n_i\in\mathbb Z_{\geq0}\ (1\leq i\leq10)\ ,\qquad
k_a-p_a^{\mathsf T}n\in\mathbb Z_{\geq0}\ ,\qquad m_a\in\mathbb R\quad(a=1,2)\ .
\end{equation}
In particular, $k_a\in\mathbb Z$. Either sign of a central charge is possible.
\end{proposition}
\begin{proof}
For each simple $\mathfrak{sl}_2$, put $\ell=\Lambda(h_i)$. Hermiticity makes $\ell$ real. Commuting $e_i$ through a power of $f_i$ gives
\begin{equation}\label{eq:norm}
  \|f_i^r v_\Lambda\|^2
  = r!\,\ell(\ell-1)\cdots(\ell-r+1)\ ,\qquad r\geq1 \ .
\end{equation}
Positivity for $r=1$ forces $\ell\geq0$. If $\ell$ were nonintegral,
$r=\lfloor\ell\rfloor+2$ would give a negative squared norm. Thus, all twelve simple labels are nonnegative integers, and the $m_a$ are real by self-adjointness. Conversely, these conditions give a dominant integral highest weight real on the chosen real Cartan. The unitarizability theorem for integrable irreducible highest-weight modules of a symmetrizable Kac--Moody algebra gives a positive-definite contravariant Hermitian form \cite[Thm.~11.7]{Kac2} (see also \cite{KP}). In the coordinates
\eqref{eq:weightcoords}, the conditions are exactly \eqref{eq:unitary}.
\end{proof}

If all twelve labels vanish, the irreducible module is the one-dimensional unitary character $\chi_m$ that is trivial on $\mathfrak g'(A_3)$ and satisfies $\chi_m(d_a)=m_a$. For a general unitary highest-weight module with degree eigenvalues $m=(m_1,m_2)$ on its highest-weight line, tensoring with $\chi_{-m}$ sets both eigenvalues to zero. This is a normalization up to
tensoring by a one-dimensional unitary character. For the full algebra with the chosen $d_a$, different pairs $m$ give inequivalent highest-weight modules at fixed Dynkin labels.
Already, $\chi_{(0,0)}$ and $\chi_{(1,0)}$ are distinguished by the action of $d_1$. If at least one Dynkin label is positive, the module is infinite-dimensional. In fact, at least one of the two affine six-cycles then has positive level on $v_\Lambda$, since its canonical central element is the sum of its six coroots. On the cyclic module for that affine derived subalgebra, the sum acts by a positive scalar. In finite dimension, the trace of each $h_i=[e_i,f_i]$ would vanish, which contradicts the trace of their
sum.

For the prism one can read off the possible central charges without solving any inequalities. The charge columns of $(v_1\ v_2)^{\mathsf T}$ at nodes $11$ and $2$ are $(1,0)$ and $(-1,0)$. At nodes $12$ and $3$ they are $(0,1)$ and $(0,-1)$. Given any
$(k_1,k_2)\in\mathbb Z^2$, assign the Dynkin labels
\begin{equation}\label{eq:lift}
 \ell_{11}=\max(k_1,0)\ ,\quad \ell_2=\max(-k_1,0)\ ,\qquad
 \ell_{12}=\max(k_2,0)\ ,\quad \ell_3=\max(-k_2,0)
\end{equation}
and take all others to be zero. This gives a dominant lift of $(k_1,k_2)$. Both $v_a$ have coordinate sum zero, so adding the same positive integer to all twelve labels gives infinitely many distinct dominant lifts, even after fixing $m_1,m_2$. For example, the label $\ell_2=1$ alone gives charges $(-1,0)$. This is where the mixed signs of the central coroots have a representation-theoretic effect. In the affine six-cycle, the central
coroot has strictly positive coefficients, and at fixed level only finitely many tuples of nonnegative Dynkin labels are possible \cite{Kac2,go}.

The same argument has a lattice formulation valid at every step $n\geq3$. Let $r=12\cdot2^{n-3}$ and identify
$Z_{\mathbb Z}=\ker (A_n^{(L)})\cap\mathbb Z^r$ with the integral centre in the simple-coroot lattice. The quotient $\mathbb Z^r/Z_{\mathbb Z}$ is free, being isomorphic to the image of an integer matrix. Thus, an integral character $Z_{\mathbb Z}\to\mathbb Z$ extends to an integral functional on the whole coroot lattice. Its coefficients need not initially be nonnegative. The constant vector is orthogonal to $Z_{\mathbb Z}$, since it is an eigenvector of the symmetric matrix $A_n^{(L)}$ of nonzero eigenvalue $-1-(n-3)L$. Adding a sufficiently large common integer to the coefficients
therefore gives a dominant integral extension with the same character. Varying that integer gives infinitely many such extensions.

\section{Conclusion}\label{sec:conclusion}
Starting from the prism construction, the order doubles at each step of the method, while the corank counts only the shifts of the seed spectrum that reach zero. This accounts for the monotone growth at $L=2$, the alternating drops at $L=1$, and the exact zeros at even $n\geq4$ for $L\geq6$. It also reconciles
the exponential growth of the nonzero coranks with their vanishing density.

In the twelve-node algebra, the same two kernel vectors determine the centre and the radical of the root lattice. The affine cycles supply actual isotropic roots, while their pairings with the opposite cycle separate them from the radical. Marking two adjacent nodes then relates the graph to the representation theory in two different ways. Each unmarked neighbour of a marked node gives a pure Serre relation of total marked degree two. The edge between the two marked nodes already gives $[e_{11},e_{12}]\ne0$ in bidegree $(1,1)$. Their mutual Serre relations first
occur in bidegrees $(2,1)$ and $(1,2)$. Determining the
irreducible constituents that survive in those two degrees is a natural next question. At the level of highest weights, the integral centre is already fully visible. Each integral character occurs for infinitely many distinct dominant weights.

\section*{Acknowledgements}
\paragraph{Funding information}
SB and VS are supported by the Interdisciplinary Thematic Institute QMat, as part of the ITI 2021--2028 program of the University of Strasbourg, CNRS and Inserm, and were supported by IdEx Unistra (ANR-10-IDEX-0002), SFRI STRAT'US project (ANR-20-SFRI-0012), and EUR QMAT (ANR-17-EURE-0024). These funding awards were made within the framework of the French Investments for the Future Program. QB is supported by IdEx funding ``Attractivity 2025'' of the University of Strasbourg. AM's work is supported by COST Action
CaLISTA CA21109 funded by COST (European Cooperation in Science and Technology).

\begingroup
\raggedright
\phantomsection
\bibliography{reffinal_v2}
\endgroup
\end{document}